\documentclass[11pt]{article}

\usepackage[latin1]{inputenc}

\usepackage{amsfonts}
\usepackage{graphics}
\usepackage{amssymb}

\usepackage[dvips]{graphicx}

\newtheorem{theorem}{Theorem}[section]

\newtheorem{corollary}[theorem]{Corollary}

\newcommand{\proof}{{\noindent \bf Proof:} }

\newcommand{\eop }{ \hfill $\Box$ }

\newcommand{\seta}{\rightarrow}

\newcommand{\gtiu}{\mathfrak{gl}(n,\R)}

\newcommand{\sltiu}{\mathfrak{sl}(n,\R)}

\newcommand{\Gl}{Gl(n,\R)}
\newcommand{\ontiu}{\mathfrak{so}(n,\R)}

\newcommand{\ptiu}{\mathfrak{p}}
\newcommand{\qtiu}{\mathfrak{q}}

\newcommand{\Prob}{\mathbf P}
\newcommand{\R}{\mathbf R}

\newcommand{\F}{\mathcal F}

\begin{document}

\begin{center}

\vspace{1cm}

 {\Large {\bf Decomposition of stochastic flows with automorphism of
subbundles component
\\[3mm]}}

\end{center}

\vspace{0.3cm}

\begin{center}
{\large { Pedro J. Catuogno}\footnote{E-mail:
pedrojc@ime.unicamp.br. 
Research partially supported by CNPq 302.704/2008-6, 480.271/2009-7
and FAPESP 07/06896-5.} \ \ \ \ \ \ \ \ \ \ \ \ \  Fabiano B. da
Silva\footnote{E-mail: fbsborges@zipmail.com.br. Research
supported by CNPQ, grant no. 142655/2005-8}

\bigskip

{ Paulo R. Ruffino}\footnote{Corresponding author, e-mail:
ruffino@ime.unicamp.br.
Research partially supported by CNPq 306.264/2009-9, 480.271/2009-7
and FAPESP 07/06896-5.}}

\vspace{0.2cm}

\textit{Departamento de Matem\'{a}tica, Universidade Estadual de Campinas, \\
13.081-970- Campinas - SP, Brazil.}

\end{center}

\begin{abstract}
We show that given a $G$-structure $P$ on a differentiable 
manifold $M$, if the group $G(M)$ of automorphisms of $P$ is large enough,
then there
exists
the quotient of an stochastic flows $\varphi_t$ by $G(M)$, in the sense 
that $\varphi_t = \xi_t \circ \rho_t$ where $\xi_t \in G(M)$, the remainder
$\rho_t$  has derivative which is vertical,
transversal to the fibres of $P$. This geometrical context generalises
previous results where $M$ is a Riemannian manifold and $\varphi_t$ is  
decomposed with an isometric component,
see \cite{Liao1} and \cite{Ruffino}, which in our context corresponds to
the particular case of an $SO(n)$-structure on $M$.
\end{abstract}

\noindent {\bf Key words:} decomposition of stochastic flows, automorphisms of
$G$-structures, infinitesimal automorphisms, stochastic exponential, symplectic
flows.

\vspace{0.3cm}
\noindent {\bf MSC2010 subject classification:} 60H10, 58J65, (53C10).

\section{Introduction}

Let $M$ be a connected differentiable manifold, with 
$\mathrm{Diff} (M)$ the group of smooth diffeomorphisms of $M$. 
Consider the following Stratonovich stochastic differential equation (sde) on
$M$:
\begin{equation}
dx_{t} \ = \ \sum_{i=0}^{k} X_{i}(x_{t})\circ dW_{t}^{i}
\label{eq1}
\end{equation}
with initial condition $x_0\in M$, where $X_0, X_1, \ldots X_k$ are smooth
vector fields on
$M$, $(W_{t}^{0})=t$, and $(W_{t}^{1}, \cdots W_{t}^{k})$ is a
Brownian motion in $\R^k$, defined over an appropriate filtered
probability space $(\Omega,\F, (\F_t )_{t\geq 0},  \Prob)$. We shall denote by 
$\varphi_t: \Omega \times M
\rightarrow M$ the stochastic flow associated to the diffusion generated by this
equation, which we shall assume to exist for all $t\geq 0$, e.g. assuming that
the derivatives of the vector 
fields of the sde are bounded.

In this article we study decompositions of the stochastic flow $\varphi_t$ such
that one of  the
components in the decomposition is a diffusion in the group of automorphisms of
a $G$-structure on $M$, i.e. a subbundle of
the principal
bundle of linear frames in each tangent space of $M$, with $G$ as the structure
group. Our main
result shows that given a $G$-structure $P$ on
a differential 
manifold $M$, if the group of automorphisms $G(M)$ is large enough, then there
exists
the quotient of the stochastic flow $\varphi_t$ by $G(M)$, in the sense 
that, for a fixed initial conditions $x_0\in M$ and a frame $u_0 \in P$,
$\varphi_t = \xi_t \circ \rho_t$ where $\xi_t \in G(M)$, the remainder
$\rho_t(x_0)=x_0$ for all $t\geq 0$ and $\rho_t$ has linearization which is
vertical, 
transversal to the fibre of $P$ at $u_0$, precisely  $\rho_{t*}(u_0)=u_0 \cdot
q_t$, where $q_t$ is a process which lives in the exponential of a
complementary subspace of the Lie algebra of $G$. We explore geometrically
interesting examples where
these complementary spaces are Lie subalgebras of $\gtiu$.

This geometrical $G$-structure context generalises
previous results of Liao \cite{Liao1} and Ruffino \cite{Ruffino} where they
decompose
$\varphi_t$ with an isometric component, which in our context corresponds to
the particular case of an $SO(n)$-structure on $M$. In those articles the
decomposition (quotient) has been used to calculate the
Lyapunov exponents and the matrix of rotation, respectively.

In the following section we shall recall the 
definitions and basic properties of $G$-structures on a differentiable manifold
$M$. The main results are presented in Section 3. In Section 4 we present a
sequence
of examples and applications. The reader will notice that the approach used here
can also be used, with few adaptations, to proof the same results for control
and nonautonomous flows. The fact that we deal with  stochastic equations will
guarantee that if the original system of equation (\ref{eq1}) has
sufficiently many
directions of diffusion such that the corresponding system  in $G(M)$ is
nondegenerate (Equation
(\ref{eq9})), then the support of $\xi_t$  will be the
connected
component of the identity of $G(M)$. We remark also that the component $\xi_t$
illustrates, for $r=1$ the 
stochastic calculus of order $r$ for diffusions in $M$ which
are automorphisms of $G$-structures of contact order $r$, cf. Akiyama
\cite{Akiyama}.

\section{Geometric set up}\label{nocoesGeometria}

In this section we introduce
 the geometrical objects involved in the general technique of
 decomposition of flows given a group of automorphisms of a
subbundle of a principal bundle, which most of interesting case are
G-structure. We refer mainly to  the classical Kobayashi
\cite{Kobayashi}
or Kobayashi and Nomizu \cite{Kobayashi-Nomizu}.

Let $M$ be a smooth connected $n$-dimensional manifold $M$. We shall denote by
$GL(M)$
the principal bundle of linear frames in each tangent space of $M$, i.e. the set
of linear isomorphisms $u:\R^{n} \rightarrow T_{x}M$ for all $x\in M$
with the natural projection $\pi:GL(M) \seta M$ where $Gl(n,\R)$ is the
structural group.

We shall consider a complete connection $\Gamma$ in $GL(M)$, i.e. such that any
segment of $\Gamma$-geodesic in $M$ defined for parameter $t$ in an interval can
be extended to all  $t\in \R$. The corresponding
$\gtiu$-value 1-form connection will be denoted by $\omega$, where $\gtiu$
denotes the Lie algebra
of $Gl(n, \R)$. For any  vector  $v \in \R^n$, there exists a standard
horizontal vector field $B(v)$ in $GL(M)$ given by the following: at a frame
$u$,
$B(v)$ is the unique horizontal vector in $T_uGL(M)$ such that $d_u \pi(B) =
u(v)$. Alternatively, if
$\theta$ is the canonical form $\theta_u : T_uGL(M) \seta
\R^n$ given by $\theta_u (Y) = u^{-1}d\pi (Y)$, then
$\theta (B(v))=v$. Let $H_1,...,H_n$ be the standard horizontal
vector fields generated by the canonical basis $e_1,...,e_n$ of
$\R^n$ and let $\{E_i^{j*}\}$ be the fundamental vector fields
corresponding to the basis $\{E_i^{j}\}$ of the Lie algebra
$\gtiu$, where the fundamental vector fields in $GL(M)$ is given, for each 
$A$ in the Lie algebra $ \gtiu$, by
\[
 A^* (u) =
\frac{d}{dt}[u\exp(tA)]_{t=0}.
\]

 These $n^2 + n$
vector fields $\{ (H_k)_u , (E_i^{j*})_u \}$ form a basis of
$T_uGL(M)$ for every $u \in GL(M)$ with $\theta (H_k)=e_k$ and $\omega
(E_i^{j*}) = E_i^{j}$.

\bigskip

A diffeomorphism $\varphi:M \seta M$ induces naturally an
automorphism $\varphi_*$ of the bundle $GL(M)$; which maps a frame
$u=(X_1,...,X_n)$ at a certain point $x \in M$ into the frame
$\varphi_*(u)=(\varphi_*X_1,...,\varphi_*X_n)$ at $\varphi(x) \in M$. The
canonical form $\theta$ is invariant 
by the pull-back $\varphi^* \theta=\theta$. In particular, we say that a
diffeomorphism $\varphi: M \seta
M$ is an affine transformation if the derivative $\varphi_*: TM \seta TM$
maps horizontal curves into horizontal curves. Equivalently, $\varphi$ is
affine if  $\varphi_*$
preserves the horizontal subspaces established by the connection $\Gamma$, or
yet, if it preserves the connection form $\varphi^*\omega=\omega$.

%% Kobay. e Nomizu pg229 %%

Let $X$ be a  smooth vector field on $M$, and let  $\eta_t$ be  its 
associated
one-parameter flow. The natural lift of $X$, denoted by $\delta X(u)$, is a
vector
field in $\mathfrak{X}(GL(M))$, the Lie algebra of vector fields in the frame
bundle, given by:
\[
\delta X(u) = \frac{d}{dt} \eta_{t*} (u) \arrowvert_{t=0}.
\]
The lift $\delta X$ can also be characterised by the following properties
simultaneously:
\begin{description}
\item a) $R_{a*}\delta X(u) = \delta X(ua) $ for every $a \in G$;

\item b) $L_{\delta X} \theta = 0$ (Lie derivative);

\item c) $d\pi (\delta X (u)) = X(\pi (u))$, for every $u \in
GL(M)$.

\end{description}

Given a complete linear connection $\Gamma$ and its connection form
$\omega$ on $M$, a vector field $X$ on $M$ is an infinitesimal
affine transformation if the associated flow $\eta_t$ are affine transformations
for all $t\in \R$. Affine infinitesimal transformations, denoted 
by $a(M)$ is a Lie algebra  isomorphic to the subalgebra of
$\omega$-preserving elements of $\mathfrak{X}(GL(M))$:
\[ 
a(GL(M))=\{ \delta X \in \mathfrak{X}(GL(M)): 
L_{\delta X}\omega=0 \}.
\]

The Lie algebra $a(M)$ has dimension at most $n^2+n$, and for any $u\in GL(M)$
the linear mapping $a(M) \seta T_uGL(M)$, given by $ X
\mapsto \delta X(u)$ is injective. When  $\dim a(M) =
n^2+n$ then $\Gamma$ is flat. See e.g.  Kobayashi and Nomizu
\cite[Chap. 3, Thm 2.3]{Kobayashi-Nomizu} among others.

 Fundamental vector fields $A^* (u)$, with $A \in \gtiu$ generate 
the vertical subspace $T_{u}^{v}GL(M)$ of $T_{u}GL(M)$. For $Y \in
T_xM$, let $\nabla_Y X= \nabla X (Y)$ be the covariant derivative
defined by the Riemannian connection. The vertical component of the
canonical lift is given precisely by $(\delta
X)^{v}(u)=\nabla X(u)$.
Hence, there is a unique matrix $[\widetilde{X}(u)] \in \gtiu$
which, acting on the right, equals the covariant derivative of $\nabla X$ acting
on the left: i.e.  $\nabla X(u) = u [\widetilde{X}(u)]$. We refer to Cordero
et al. \cite{Cordero}.

%\section{Framework}

\subsection{$G$-structures}\label{Gestrutura}

A $G$-structure on $M$ is a reduction of the frame bundle
 $GL(M)$ to a subbundle $P$, with structure Lie group
 $G \subset
Gl(n,\R)$. Given a closed subgroup $G$, the existence 
of a $G$-structure
depends intrinsically on the topology of $M$: More precisely,
 there exists a one to one correspondence between $G$-structures
 and  cross sections of the associated fibre bundle  $GL(M)/G$, see
\cite[Chap. I, Prop. 5.6]{Kobayashi-Nomizu}. 
A list of interesting examples includes the following: There 
exist $Gl(n,\R)^+$-structures if and only if $M$ is orientable; there 
exist $ Sl(n, \R)$-structures if and only if $M$ has a volume
 form; paracompactness guarantees the existence of
$O(n, \R)$-structures via Riemannian metrics; finally note that
 a manifold $M$ is parallelizable if  and only if there exist 
a $\{1\}$-structure. For these properties and further example,
 we refer to \cite{Kobayashi}, \cite{Kobayashi-Nomizu}, the classical
Sternberg 
\cite{Sternberg}, among others.

\bigskip

 Given a $(k,l)$-tensor $K$ over the Euclidean space
$\R^n$, let $G$ be the group of linear transformations in
$\R^n$ which is $K$-invariant, i.e. for $(v_1,...,v_k, f_1, \dots , f_l) \ \in
(\R^n)^k\times (\R^{n*})^l$, we have that $g \in G$ if
\[
K(v_1,...,v_k, f_1, \dots , f_l)= K(gv_1,...,gv_k, g^*f_1, \dots , g^*f_l).
\]
We say that a corresponding $G$-structure over $M$ is induced by the tensor $K$.
A such $G$-structure extends naturally the tensor $K$  to a tensor field
$k$ on $M$ defined by: for each $x \in M$, given a
linear isomorphism $u:\R^n \seta T_xM \in P$, for
$w_1,...,w_k \ \in T_xM$ and $z_1,...,z_l \ \in T_xM^*$, we assign
\[
k_x(w_1,...,w_k, z_1,...,z_l ) =  K(u^{-1}w_1,...,u^{-1}w_k,
u^{-1*}z_1,..., u^{-1*}z_l).
\]
The invariance of $K$ by $G$ guarantees that the definition above is independent
of the choice of $u$.

\bigskip

Let $\varphi:M \seta M$ be a smooth diffeomorphism, given a
$G$-structure $P$, if
$\varphi_*$ maps
$P$ into itself, we call $\varphi$ an automorphism of the $G$-structure $P$.
A vector field $X$ on $M$ is called an infinitesimal automorphism
of a $G$-structure $P$ if it generates a local 1-parameter group
of automorphisms of $P$. A vector field $X$ on
$M$ is an infinitesimal automorphism if and only if $L_X k=0$,
where $k$ is the tensor field associated to the $G$-structure $P$.

To illustrate, we recall that $Gl (n, \R)^+$-structures has the group of
automorphisms given
 by diffeomorphisms which preserve orientation, and the
 infinitesimal automorphisms are differentiable vector
 fields. For $ Sl(n, \R)$-structures, the group of 
automorphisms is given by diffeomorphisms which preserve volume,
 and the infinitesimal automorphisms are vector fields with
 vanishing divergent. For $O(n, \R)$-structures, the group 
of automorphisms is given by isometries and the infinitesimal
 automorphisms are Killing vector fields. For a $\{1\}$-structure,
 the group of automorphisms degenerates to identity and the
 infinitesimal automorphisms degenerates to the zero vector field.

%%%%%%%%%%%%%%%%%%%%%%%%%%%%%%%%%%%%%%%%%%%%%%%%%%%%%%%%%%%%%%%%%%%%%%%%%%%%%%%%
%%%%%%%%%%

\section{Main Results}

We shall consider a differentiable manifold $M$ 
 which admits a $G$-structure $P$
induced by a tensor field $k$, where $G$ is a subgroup of $\Gl$. 
Let $\pi_0:P \seta M$ be the restriction of $\pi:
GL(M) \seta M$. We shall fix an initial condition $x_0 \in M$ and an initial
frame $u_0 \in P$ in the tangent space at $x_0$, i.e. 
$\pi_0 (u_0)=x_0$.

We assume that $M$ is endowed with a complete
connection $\nabla$. Let $ g(M)= \{ X \in a(M): L_X k=0 \} $, since the Lie
derivative $L_{[X,Y]}= [L_x, L_Y]$ (see e.g. \cite[Prop 3.4
p.32]{Kobayashi-Nomizu}), we have that $g(M)$ is a Lie algebra
of affine
 vector fields which preserve the tensor $k$. Given a vector field $X$ in
$g(M)$, the
associated flow $\eta_t$ is a family of automorphisms of $P$. Hence
 $\eta_{t*}(u_0) \in P$
consequently, $\delta X(u_0) \in T_{u_0}P $. We shall denote by
$\gamma$ the restriction of the canonical lift $\delta $ to
vector fields in $g(M)$:  
\begin{eqnarray}
\gamma : g(M) & \seta & T_{u_0}P \nonumber \\
      X    & \mapsto & \delta X(u_0) \label{i1}.
\end{eqnarray}
It is injective since
 it is the restriction of the linear and injective
mapping $a(M) \seta T_{u_0}GL(M)$, $X \mapsto \delta X(u_0)$.

The finite dimensional Lie group generated by $g(M)$ will be denoted by
$G(M)$. So, if $\xi_t$ is a process in $G(M)$ then $\xi_t$ is an affine
transformation which preserves the tensor field $k$. 

We shall use freely the Lie group terminology to deal with flows and vector
fields, e.g. given $\xi \in G(M)$ and $Y \in g(M)$, $\xi Y$ and $Y\xi$ are
the tangent vectors at $T_{\xi}G(M)$ obtained by left and right
translations of $Y$ respectively. Clearly, $\xi \mapsto \xi Y $ are
left invariant vector field on $G(M)$. Analogously for the right
translation $Y \xi$.

Given initial conditions $x_0\in M$ and $u_0 \in \pi^{-1}(x_0) \subset P$
 we shall fix a projection $p:T_{u_0}GL(M) \rightarrow T_{u_0}P$ such that
 $p ( HT_{u_0}GL(M)) = HT_{u_0} P$ and $p (VT_{u_0}GL(M)) = VT_{u_0} P$. By
 the linear dependence of the fundamental vector fields with the elements in
the Lie algebra $\gtiu$,
 such a projection $p$
 is equivalent, and will be identified with the same notation, of a
projection in the
 corresponding Lie algebras $p: \gtiu \rightarrow \ptiu$, where 
 $\ptiu$ is the Lie algebra of $G$.

We shall assume the following hypotheses:

\begin{description}
\item[(H1)] (Existence of a $G$-structure) The manifold $M$ admits a
$G$-structure $P$. Given initial conditions $x_0\in M$ and $u_0 \in
\pi^{-1}(x_0) \subset P$, we shall fix a projection $p:T_{u_0}GL(M)
\rightarrow
T_{u_0}P$ such that
 $p ( HT_uGL(M)) = HT_u P$ and $p (VT_uGL(M)) = VT_u P$.

% i.e.:
%\[ 
%           T_{u_0}GL(M) \cong T_{u_0}P \oplus u_0 \qtiu
%\] 
%where $u_0 \qtiu=\{u_0A: A \in \qtiu \}$.
%
%For $Y \in
%T_{u_0}GL(M)$, let $Y_p$ be the $T_{u_0}P$-component of $Y$ in the
%above decomposition.

\item[(H2)] (The Lie algebra $g(M)$ is big enough) The projections
$p[\delta  Ad (\xi) X_j (u_0)]$ are in the image $ Im(\gamma)$, where $\gamma$
is given in equation (\ref{i1}), for all vector
fields
$X_j$,
$j=0,1,...,k$ of the sde $(\ref{eq1})$ and for any automorphism $\xi$ in $G(M)$.

\end{description}

\noindent The geometric and dynamical meaning of different choices of 
projection $p$ in (H1) above will be clear in the corollaries and examples
after  the
main result (Theorem \ref{geral}) below. This theorem 
 generalise to automorphisms of $G$-structures the 
factorizations
 presented in Liao \cite{Liao1}, Ruffino \cite{Ruffino},
 Colonius and Ruffino \cite{Colonius-Ruffino} for stochastic and control flows.

\begin{theorem}[Decompositions in automorphisms of
G-structure]\label{geral}
Assume \\ the conditions (H1) and (H2) above. Then the stochastic flow
$\varphi_t$
of
equation $(\ref{eq1})$
decomposes  as $\varphi_t= \xi_t \circ \rho_t$,
where $\xi_t$ is a diffusion in the group $G(M)$ of automorphisms of $P$, the
remainder $\rho_t$ is a process in
$\mathrm{Diff} (M)$ such that $\rho_t(x_0)=x_0$ and
$\rho_{t*}(u_0)= u_0 \cdot q_t $, where $q_t$ is a process in $ \langle
\exp(\ker p) \rangle$, the Lie
subgroup generated by $\exp(\ker p)$.
\end{theorem}

\proof

Let $\xi_{t}$ be the solution of the following sde in the group
$G(M)$ with initial
condition 
$\xi_{0} = id$:
\begin{equation} \label{eq9}
d\xi_{t} = \sum_{j=0}^{k} \left( L_{\xi_{t}}\right)_*\ [Ad(\xi^{-1}_t) X_{j}]^p
\circ
dW_{t}^{j} \ ,
\end{equation}
where $[Ad(\xi^{-1}_t)X_{j}]^p$ is the unique vector field in
$g(M)$ such that $ \delta [Ad(\xi^{-1}_t)X_{j}]^p(u_0) =
p (\delta Ad(\xi^{-1}_t)X_{j}(u_0))$. 
Note that the horizontal
components of these two vectors coincide, i.e. in $T_{x_0}M$ we have that
$[Ad(\xi^{-1}_t)X_{j}]^p(x_0)=Ad(\xi^{-1}_t)X_{j}(x_0)$.
The process $\xi_t$ is an stochastic exponential in the sense of Hakim-Dowek and
 L\'epingle \cite{H-DL} or Catuogno and Ruffino \cite{Catuogno-Ruffino}.
Since the process is the stochastic exponential of a process in the Lie algebra
$g(M)$, then $\xi_{t}$ is a diffusion in the
group $G(M)$ of
automorphisms of
$P$.

For the remainder $\rho_t = \xi_t^{-1} \varphi_t$, by the Itô formula in the
group, as in Kunita \cite[pp. 208-209]{Kunita} and  the
fact that
\begin{equation} 
d\xi_{t}^{-1} = - \sum_{j=0}^{k} \ [Ad(\xi^{-1}_t)X_{j}]^{p}
\ \xi_{t}^{-1} \circ dW_{t}^{j} \ ,
\end{equation}
we have that 
\begin{equation} 
d\rho_t = \sum_{j=0}^{k} \{ Ad(\xi^{-1}_t)X_{j} -[Ad(\xi^{-1}_t)X_{j}]^{p}
 \} (\rho_t) \circ dW_{t}^{j}.
\end{equation}
Since at $x_0$, the vector field
of the
equation above
vanishes, the 
derivative process $\rho_{t*}
u_0$ starting at $u_0$, has no horizontal component and
satisfies the linear sde in
$T_{u_0}GL(M)$:
\begin{equation} 
d \rho_{t*} = \sum_{j=0}^{k} \delta \{ (
\xi_{t*}^{-1}X_{j}) - [\xi_{t*}^{-1}(X_{j})]^{p}
 \} \rho_{t*}  \circ dW_{t}^{j}.
\end{equation}

The derivative $\rho_{t*}$ at $u_0$ in the fibre $\pi^{-1}(x_0)$ acts on
the left
as a linear
 transformation which has no horizontal component neither vertical component
along the action of $\ptiu$. Hence it
has a vertical component in $T_{u_0}GL(M)$ along a process which is in the
Lie group generated by the kernel of the projection $p$, i.e. $\rho_{t*} u_0 =
u_0 \cdot q_t$, where $q_t$ is in $ \langle \exp(\ker p) \rangle \subseteq
\Gl(n)$.

\eop
  
In other words, the theorem above says that there exists a right quotient of a
stochastic flow in
the differentiable manifold $M$ with
respect to the group of automorphisms of a $G$-structure if the
adjoint by automorphisms of each vector
field involved in the sde is equal to a certain infinitesimal transformations
at $x_0$ and both have the same canonical lift to $u_0$. The
``remainder'' of the quotient is the process $\rho_t$.

\bigskip

If the Lie algebra $\ptiu$ has a complementary Lie subalgebra
$\qtiu$, i.e. $\gtiu= \ptiu
\oplus \qtiu$, there is a natural choice for the projection $p$ in hypothesis
(H1). In this case, the derivative of the  remainder $\rho_t$ becomes uniquely
characterised at
the initial condition.

\begin{corollary}\label{decompoe gln}
If the Lie algebra $\ptiu $ of the Lie group $G$ associated to the $G$-structure
$P$ has a complementary Lie algebra $\qtiu$,  then 
there exists a decomposition of the stochastic flow $\varphi_t$ of
equation $(\ref{eq1})$ as $\varphi_t= \xi_t \circ \rho_t$,
where $\xi_t$ is a diffusion in the group $G(M)$ of automorphisms of $P$, 
$\rho_t$ is a process in
$\mathrm{Diff} (M)$ such that $\rho_t(x_0)=x_0$ and
$\rho_{t*}(u_0)=  u_0 \cdot q_t $, where
 $q_t$ is a
process in the Lie
group $ \langle \exp\qtiu \rangle $ generated by the Lie algebra $\qtiu$.

The processes $q_t$ which satisfies the property above is unique. 
\end{corollary}

\proof Just take the projection  $p:
\gtiu \rightarrow \ptiu$ of hypothesis (H1) along the subspace $\qtiu$ and
apply Theorem \ref{geral}. For the uniqueness of the process $q_t$ in the group
$ \langle \exp\qtiu \rangle $, let $ \xi_t \circ \rho_t$ and $\widetilde{\xi}_t
\circ
\widetilde{\rho}_t$ be two distinct decomposition of $\varphi_t$ with
$\widetilde{\rho}_{t*}(u_0)= u_0 \widetilde{q}_t$ and
$\rho_{t*}(u_0)= u_0 q_t$. We have that 
$\xi_t^{-1} \circ \widetilde{\xi}_t$ is an automorphism of the
$G$-structure
$P$ which fix the point $x_0$ for all $t\in \R$. Hence, the derivative 
$(\xi_t^{-1} \circ \widetilde{\xi}_t)_* (u_0)$ is a vertical translation
in the fibre $\pi^{-1}(x_0) \subset P$. On the other hand, since
the action of $\rho_{t*}$ is equivariant:
\begin{eqnarray*}
(\xi_t^{-1} \circ
\widetilde{\xi}_t)_* (u_0) &=& (\rho_t \widetilde{\rho}_t^{-1})_*(u_0)\\
    &=& u_0 q_t \widetilde{q}_t^{-1}.
\end{eqnarray*}
The unique element which is a vertical translation of $u_0$ in $P$ and has the
form above with $q_t \widetilde{q}_t^{-1}\in \langle \exp \qtiu \rangle$
is $u_0$
itself, it follows that $ \widetilde{q}_t = q_t$.

\eop

The uniqueness stated above with local properties does not imply that the
decomposition $\varphi_t= \xi_t
\circ \rho_t$ is unique with these properties  if the group of
diffeomorphisms of $M$ which preserves the tensor associated with the
$G$-structure is large enough. One may find a, say deterministic curve of
tensor preserving diffeomorphisms $\eta_t$, with  $\eta_0=Id$ and which
restricted to an open neighbourhood of $x_0$ is the identity for all $t\geq 0$.
In this case $\widetilde{\xi}_t = \xi_t \eta_t$ and $
\widetilde{\rho}_t = \eta_t^{-1} \rho_t$ also satisfies the local conditions
stated in the corollary, although in this case $\widetilde{\xi}_t$ may no
longer be an affine transformation. For example, in the infinite
dimensional group  of
diffeomorphisms which preserve volume we do not have uniqueness. In
the isometry group the uniqueness holds as stated in \cite{Liao1} and
\cite{Ruffino}.

\bigskip

The  diffusion of $G$-automorphism $\xi_t$ in the main theorem induces a Markov
process in $M$ which is not necessarily time homogeneous, see e.g. Liao
\cite{Liao3}. The next result shows an alternative decomposition
$\varphi_t=\xi_t \rho_t$ where $\xi_t$ is a flow (time homogeneous) in $M$
itself, instead of in $G(M)$:

\begin{theorem}
 Under conditions (H1) and (H2) above, the stochastic flow $\varphi_t$ of
equation $(\ref{eq1})$
factorizes as $\varphi_t= \xi_t \circ \rho_t$,
where $\xi_t$ is an stochastic flow in $M$ which preserves the $G$-structure.
The remainder component $\rho_t$ is a process in
$\mathrm{Diff} (M)$ such that the vector fields of its sde vanish at
$\xi_t^{-1}(x_0)$.
\end{theorem}

\proof
 Define $\xi_t$ as the solution of the following right invariant stochastic
equation. 
\begin{equation} \label{eq100}
d\xi_{t} = \sum_{j=0}^{k} [ X_{j}]^p \xi_t
\circ
dW_{t}^{j} \ ,
\end{equation}
where, as in the proof of the theorem,  $[X_{j}]^p$ is the unique
vector field in
$g(M)$ such that $ \delta [X_{j}]^p(u_0) =
p (\delta X_{j}(u_0))$. 
Again, note that the horizontal
components of these two vectors coincide, i.e. in $T_{x_0}M$ we have that
$[X_{j}]^p(x_0)=X_{j}(x_0)$.

For the remainder $\rho_t = \xi_t^{-1} \varphi_t$, again by the Itô formula we
have that 
\begin{equation} 
d\rho_t = \sum_{j=0}^{k}  Ad(\xi^{-1}_t) \left( X_{j} -[X_{j}]^{p} \right)\
\rho_t
 \  \circ dW_{t}^{j}.
\end{equation}
Last property of the statement follows directly from this equation.

\eop

%The vertical component of the linearization $T_{\xi_t^{-1}(x_0)}\rho_t:
%T_{\xi_t^{-1}(x_0)}M \rightarrow T_{\xi_t^{-1}(x_0)}M$ can be described
%briefly by: if  $v_t = T_{\xi_t^{-1}(x_0)}\rho_t (v_0)$, then, the covariant
%equation for $v_t$ is given by:
%\[
% Dv_t = \sum_{j=0}^{k}  \nabla_{x_0} \xi^{-1}_t \left( X_{j} -[X_{j}]^{p}
%\right)(v_t)\  \circ dW_{t}^{j}.
%\]
%which shows that in this case, with the flow in the manifold instead of in the
%group of diffeomorphisms, as one expects, we miss the derivative
%properties of the remainder in terms of exponential of the kernel of the
%projection $p$ in the hypothesis (H1). 

 The same kind of decomposition (quotient by the group) described above can 
also be
 performed on the left hand side instead of the right quotient as
considered in  Theorem
\ref{geral}. In this case the fixed point of the
remainder
$\tilde{\rho}_t$ turns out to be the moving point $\xi_t(x_0)$ for all $t\geq
0$: 

\begin{corollary}
With the same hypotheses of Theorem \ref{geral} we have the following left
remainder decomposition
$\varphi_t=\tilde{\rho}_t \circ \xi_t$ where $\xi_t$ is a diffusion in the group
$G(M)$ of automorphisms of the
$G$-structure, the remainder component $\tilde{\rho}_t$ is a process in
$\mathrm{Diff} (M)$ such that $\tilde{\rho}_t(\xi_t(x_0))=\xi_t(x_0)$. 
\end{corollary}

\proof Take the diffusion $\xi_t$ in Theorem \ref{geral}, i.e. equation
($\ref{eq9}$). From the
theorem we have that $\varphi_t(x_0)=\xi_t(x_0)$. Hence, for any time  $0 \geq t
$ the left hand side
remainder considered here $\tilde{\rho}_t =\varphi_t \xi^{-1}_t$,  is a
random diffeomorphism which fixes $\xi_t(x_0)$ a.s., moreover it
satisfies the non-autonomous equation.
\[
 d\tilde{\rho}_t = \sum_{j=0}^{k} \{ X_{j}
-Ad(\varphi_t)[Ad(\xi^{-1}_t)X_{j}]^{p}
 \} (\rho_t) \circ dW_{t}^{j}.
\]

\eop

\noindent \textit{Remark:} It might be possible to obtain, locally, analogous
results of this section with a 
linear connection by using local affine
transformations. This approach 
would demand to
reconstruct locally the theory of (global) affine
transformations (as in e.g. Kobayashi and Nomizu
\cite[pp.234-235]{Kobayashi-Nomizu}).
To prevent the reader from further geometrical technicalities, here we have assumed
completeness of the connection.

\section{Examples}

 We recall initially two interesting particular cases of decomposition in the
literature where
the hypothesis of Corollary \ref{decompoe gln} holds due to global geometrical
properties of the manifold. In the first case, the $G$-structure preserves
the
metric tensor, i.e. $g(M)$ is the algebra of Killing vector fields:

\begin{theorem}[\cite{Liao1}, \cite{Ruffino}] \label{decomp isom}
 If $M$ is a simply connected Riemannian manifold with constant curvature, then
every stochastic flow of the Stratonovich equation (\ref{eq1}) satisfies the
hypotheses of Corollary \ref{decompoe gln}, hence, for each initial condition
$x_0 \in M$ and $u_0$ an orthonormal frame in $T_{x_0}M$, the flow $\varphi_t$
admits the decomposition: 
\[
 \varphi_t= \xi_t \circ \rho_t,
\]
where $\xi_t$ are isometries of $M$, $\rho_t (x_0)=x_0$ and the
linearization $\rho_{t*} (u_0)= u_0 q_t$, where $q_t$ is a process in the group
of upper triangular matrices.
\end{theorem}
\proof  In this case the dimension
of $g(M)$ is maximal $n(n+1)/2$, (see e.g.
\cite{Kobayashi}) and $\gtiu = \ontiu \oplus \qtiu $, with $\qtiu$ the Lie
algebra of upper triangular matrices.

\eop

This decomposition is used to study in the same context the radial
and
angular asymptotic
behaviour: The  Lyapunov exponents comes from the remainder
$\rho_t$,
\cite{Liao1}; and the matrix of rotation comes from the isometries $\xi_t$,
\cite{Ruffino}.

The second particular case refers to
decompositions with affine transformation component, i.e. $\xi_t$ is a diffusion
in the group of diffeomorphisms which preserves the connection.

\begin{theorem}[\cite{Ruffino}] \label{decomp afim}
 In Euclidean spaces, given the stochastic flow of the Stratonovich equation
(\ref{eq1}) for each initial condition
$x_0 \in M$ and $u_0$ a frame in $T_{x_0}M$, then
\[
 \varphi_t= \xi_t \circ \rho_t,
\]
where $\xi_t$ are affine transformations, $\rho_t (x_0)=x_0$ and the
linearization $\rho_{t*} = Id$. 
\end{theorem}
\proof  It is a direct consequence of Corollary \ref{decompoe gln}. Again, the
dimension of $g(M)$ is maximal $n(n+1)$, (see e.g.
\cite{Kobayashi-Nomizu}) and $\ptiu = \gtiu$, hence it degenerates the
derivative of the remainder $\rho_t$. 

\eop

 The theorem above holds in a more general geometrical context: essentially the
group
of affine transformations has to be large enough (depending also on the sde
considered). Again, for sufficiently large groups of affine transformations of
compact manifolds, Theorem \ref{decomp afim} implies Theorem \ref{decomp isom},
since in this case, affine transformations are isometries (see e.g. Kobayashi
\cite[Cor.2.4]{Kobayashi}).

\subsection{Decomposition of symplectic flows}

%We denote by $Sp(n, \R)= \{ A
%\in Gl(2n, \R): A^tJA = J \}$ the symplectic group, where
%$$J=\left(
%\begin{array}{cc}
%  0 & -I  \\
%  I & 0 \\
%\end{array}
%\right).$$
%Its associated Lie algebra is given by $\mathfrak{sp}(n,\mathbb{R})= \{ A
%\in \mathfrak{gl}(2n,\mathbb{R}): A^tJ+ JA = 0 \}$. Hence,
%$\mathfrak{sp}(n,\mathbb{R})$
%consists of matrix of the form
%$$A=\left(%
%\begin{array}{cc}
%  A_1 & A_2  \\
%  A_3 & -A_1^t \\
%\end{array}%
%\right), \ \textrm{with} \  A_2^t= A_2 \  \textrm{and} \ A_3^t= A_3.$$

Let $(M, \omega)$ be a symplectic manifold, where $\omega$ here
denotes a closed nondegenerate 2-form in $M$. In this example we consider  a
stochastic symplectic flow $\varphi_t$ in $M$
associated with the sde
\begin{equation} \label{eq1s}
dx_{t} \ = \ \sum_{i=0}^{k} X_{i}(x_{t})\circ dW_{t}^{i}
\end{equation}
where $X_0, X_1, \ldots X_k$ are smooth symplectic vector fields on
$M$. The flow $\varphi_t$ is a symplectic
transformation of $M$, i.e., $\varphi_t^* \omega = \omega$ almost surely (see
e.g. Kunita \cite{Kunita}).

Consider the following decomposition of the symplectic Lie algebra 
$\mathfrak{sp}(n,\R) = \ptiu \oplus \qtiu$ where
$$\ptiu= \bigg\{ \left(
\begin{array}{cc}
  A & -S  \\
  S & A \\
\end{array}
\right):  \  A^t= -A \  \textrm{and} \ S^t= S \bigg\}.$$
And 
$$\qtiu= \bigg\{ \left(
\begin{array}{cc}
  \Delta & S  \\
  0 & -\Delta^t \\
\end{array}
\right):  \  S^t= S \mbox{ and } \Delta \mbox{ is upper
triangular}\bigg\}.$$

Given an element $A \in \mathfrak{sp}(n,\R)$, we write
$A= A_{\mathfrak{p}}+ A_\mathfrak{q} $ with
$A_{\mathfrak{p}} \in \mathfrak{p}$ and $A_{\mathfrak{q}} \in \mathfrak{q}$.
Note that the Lie algebra $\mathfrak{p}$ is a Lie subalgebra of
skew-symmetric matrices. Therefore, the subgroup generated 
$\langle \exp \ptiu \rangle \subset Sp(2n, \R)\cap SO(2n, \R)$. We
assume that
$M$ admit
a $\langle \exp \ptiu \rangle$-structure $P$, with $P$ been a differentiable
subbundle of the symplectic $Sp(n)$-structure of $GL(M)$.

We consider the projection $p: \mathfrak{sp}(n,\R) \rightarrow \ptiu $
along the subspace $\qtiu$. We have a direct consequence of Corollary
\ref{decompoe gln}:

\begin{corollary} \label{corol sympl}
The symplectic stochastic flow $\varphi_{t}$ has a decomposition
$\varphi_{t} = \xi_{t} \circ \rho_{t}$, where $\xi_{t}$ is a diffusion
in the group of isometries, $\rho_{t}(x_0) = x_0$ and
$\rho_{t*}(u_0) = u_0q_{t}$, for some process $q_{t} $ in the subgroup $\exp
(\qtiu)$ of the symplectic group.
\end{corollary}

This decomposition with the choice of subalgebras $\ptiu$ and
$\qtiu$ as above is not the same decomposition with isometric component of
Theorem \ref{decomp isom}. To illustrate, consider this simple example:

$$ \dot{x_t} = u_1(t) A (x_t)+u_2(t)B(x_t) $$
with non commutative 
$$A=  \left(
\begin{array}{cc}
\left(
\begin{array}{cc}
  0 & 1  \\
  0 & 0 \\
\end{array}
\right)
& 0  \\
  0 &
\left(
\begin{array}{cc}
  0 & 0  \\
  -1 & 0 \\
\end{array}
\right)
\\
\end{array}
\right)  \ \textrm{ and } \ \
B= \left(
\begin{array}{cc}
\left(
\begin{array}{cc}
  0 & -1  \\
  1 & 0 \\
\end{array}
\right)
& 0  \\
  0 &
\left(
\begin{array}{cc}
  0 & -1  \\
  1 & 0 \\
\end{array}
\right)
\\
\end{array}
\right)
$$
and the bounded measurable control functions given by the indicator functions
of intervals $u_1(t)= 1_{[0,1]}(t)$, $u_2(t)=1_{[1,\infty]}(t)$. The initial
conditions are $x_0$ the origin and $u_0$ the canonical basis of $\R^{2n}$.

The solution of the control equation above is the symplectic diffeomorphism
$\varphi_t$, which has the following decomposition: For $0\leq t\leq 1$
\begin{equation} \label{dec sympl1}
 \varphi_t = Id \circ 
\left(
\begin{array}{cc}
\left(
\begin{array}{cc}
  1 & t  \\
  0 & 1 \\
\end{array}
\right)
& 0  \\
  0 &
\left(
\begin{array}{cc}
  1 & 0  \\
  -t & 1 \\
\end{array}
\right)
\\
\end{array}
\right)
\end{equation}
and for $ t\geq 1$,
 
\begin{equation}   \label{dec sympl2}
\varphi_t= \left(
\begin{array}{cc}
\mathrm{Rot} (t-1)
& 0  \\
  0 &
\mathrm{Rot} (t-1)
\\
\end{array}
\right)
\left(
\begin{array}{cc}
\left(
\begin{array}{cc}
  1 & 1  \\
  0 & 1 \\
\end{array}
\right)
& 0  \\
  0 &
\left(
\begin{array}{cc}
  1 & 0  \\
  -1 & 1 \\
\end{array}
\right)
\\
\end{array}
\right),
\end{equation}
where the 2-dimensional rotation 
\[ \mathrm{Rot (s)}:=
\left(
\begin{array}{cc}
  \cos s & -\sin s  \\
  \sin s  & \cos s \\
\end{array}
\right),
\]
for $s\in \R$. By commutativity of the initial frame $u_0$ with each component,
one sees
that the first matrices in the product of equations (\ref{dec sympl1}) and
(\ref{dec sympl2}) are the $\xi_t$ component of the decomposition of
Corollary
\ref{corol sympl} and the
second matrices are the remainder $\rho_t$. This decomposition differs from the
decomposition of  Theorem
\ref{decomp isom} since in the example here the derivative of the remainder in
$u_0$ is not described by the right action of an upper triangular matrix.

\subsection{Volume preserving component}

\label{Sl-structure}

in this example we give detailed calculation for the decomposition with a
volume preserving component. Consider initially that our sde of equation
(\ref{eq1}) is in the
Euclidean space $\R^n$. We shall consider the canonical volume tensor $v=dx_1
\wedge \ldots dx_n$ and the corresponding infinite dimensional group of
volume preserving diffeomorphisms. The theory of the
previous section applies here to find a decomposition of the solution flow 
 $\varphi_t$ as $\xi_t \circ \rho_t$, where the component $\xi_t$ is
in the intersection of affine transformations and the volume preserving group of
diffeomorphisms. In our terminology, it corresponds precisely to work with
the $Sl(n, \R)$-structure, which here trivializes  as $P=\R^n \times
\mathrm{Sl}(n, \R)$.

 Given an initial condition $x_0$ and an initial frame $u_0 \in P$, we consider
the following basis for the tangent space $T_{u_0}P$:
$$
\{ Y_1,...,Y_n,Y_{n+1},...,Y_{(n^2+n-1)} \},
$$ 
where, for $j=1, \ldots, n$, the horizontal elements are $Y_j=e_j$,
the canonical basis and for $j=(n+1), \ldots , (n^2+n-1)$, $Y_j=A_j$, with
$(A_j)$ a basis of $\sltiu$.

%Consequently, we have that
%$\textrm{\textrm{div}} Y_j=0$, for $j=1,...,n,...,n^2+n-1$.

%Besides, by choose of Lie algebra $g(\mathbb{R}^n)$, its easy
%checks that, for $u_0 \in P$,
%\begin{eqnarray}
%\gamma : g(\mathbb{R}^n) & \seta & T_{u_0}P \nonumber \\
%      Y    & \mapsto & \delta Y(u_0) \label{i4}
%\end{eqnarray}
%is a bijective linear map and that linear combination of vector
%fields $$\{ \delta Y_1(u_0),...,\delta Y_n(u_0),\delta
%Y_{n+1}(u_0),...,\delta Y_{n^2+n-1}(u_0)  \}$$ generates $T_{u_0}P$. More
%precisely, $\{ \delta Y_1(u_0),...,\delta Y_n(u_0) \}$ generates the
%horizontal subspace and $\{ \delta Y_{n+1}(u_0),...,\delta
%Y_{n^2+n-1}(u_0) \}$ the vertical subspace of $T_{u_0}P$. Hence, the
%hypothesis (H2) of Theorem \ref{geral} is automatically satisfies.

For $A \in \gtiu$ there is a unique decomposition $\ptiu \oplus \qtiu$ given by
$$
 A =(A - \frac{\mathrm{tr}(A)}{n}Id) + \frac{\mathrm{tr}(A)}{n}Id=
A_{\ptiu}+A_{\qtiu}
$$
where $\ptiu =\sltiu$ and $\qtiu\simeq\R$ are Lie algebras of matrices.
Hence Corollary \ref{decompoe gln} applies and we get:

\begin{corollary} \label{pres vol}
The stochastic flow $\varphi_{t}$ has a unique decomposition
$\varphi_{t} = \xi_{t} \circ \rho_{t}$, where $\xi_{t}$ is a
diffusion process in the intersection of affine transformations and
the group of volume preserving diffeomorphisms, $\rho_{t}(x_0) = x_0$ and
$\rho_{t*}(u_0) = q_{t} u_0$, for some process $q_{t} \in \R $.
The process $q_t$ carries the information of the sum of the Lyapunov exponents
of the flow $\varphi_t$:
\[
 \lim_{t\rightarrow \infty} \frac{1}{t} \log q_t = \sum_{i=1}^{n} \lambda_i
\]
where $\lambda_i$ are the $n$ Lyapunov exponents (with possible repetition).
\end{corollary}

\proof The first part follows directly from Corollary \ref{decompoe gln} and
the commutativity of the elements of
$ \qtiu$. The second statement follows because 
\[
 \sum_{i=1}^{n} \lambda_i = \lim_{t\rightarrow \infty} \frac{1}{t} \log \|
\varphi_{t*}\|,
\]
and $q_t=\| \varphi_{t*}\|$. See e.g. Arnold \cite{Arnold}, Baxendale
\cite{Baxendale}, Elworthy \cite{Elworthy}, Liao \cite{Liao2} among others.

\eop

More generally, for a Riemannian manifold $M$ that admits an
$Sl(n,\R)$-structure $P$, we have similarly the following decomposition
of the vertical
component of the natural lift of $X$:
$$ 
\nabla X=(\nabla X - \frac{ \mathrm{div}X}{n} Id)+ \frac{\mathrm{div} X}{n} Id
, 
$$
So, the same decomposition and Lyapunov
property of Corollary \ref{pres vol} holds. 

%\begin{remark}
%If we take $P=\mathbb{R}^n \times SO(n, \R)$, the $SO(n,
%\R)$-structure on $\mathbb{R}^n$, and we construct the Lie
%subalgebra $g(\mathbb{R}^n)$ in the analogous manner way to case
%$Sl(n, \R)$-structure above, the Liao decomposition follows
%naturally, because the hypothesis (H2) is automatically satisfies.
%
%\end{remark}
%%%%%%%%%%%%%

\subsection{Cascade Decompositions}

%If the Lie subalgebra $g(M)$ is $a(M)$ itself, its decomposition
%\rho_t$, with $\xi_t$ being a diffusion on affine transformations
%group of $M$, and $\rho_t$ has the following characterization:
%$\rho_t(x_0)=x_0$ and $\rho_{t*}=Id_{T_{x_0}M}, \ \forall t \geq
%0$. This particular case has been presented in Ruffino
%\cite{Ruffino}. Through this fact and by using Theorem
%\ref{geral}, we can obtain a cascade of factorization for the
%stochastic flow $\phi_t$. More precisely, we need the following
%conditions for cascade:

In the special geometrical conditions where there exists a sequence of of
subbundles $P^1 \subset P^2  \ldots \subset P^n \subset GL(M)$ which are all
$G$-structure, repeating conveniently the decomposition technique described in
the
previous section  allows a cascade (geometrical filtration) of
decompositions.
We denote, as before, by $G^i(M)$ and $g^i(M)$ the subgroup of automorphisms of
$P^i$ and the Lie algebra of infinitesimal automorphisms,
respectively for $i=1,..,n$.  The Lie
algebra of infinitesimal automorphisms is a flag of subalgebras of
vector fields
$$
g^1(M) \subset ... \subset g^{n-1}(M) \subset g^n(M) 
$$
where $g^i(M)=\{ X\in a(M): L_X k^i=0 \}$, $i=1,..,n$, for some
tensor field $k^i$. 

We shall assume that, for  $i=1,2, \ldots, n$:

\begin{description}
    \item[(C1)] There exists a sequence of $G$-structure $P^1 \subset
P^2  \ldots  \subset P ^{n} \subset GL(M)$ such that for an initial condition
$x_0\in M$ and $u_0 \in P^1$, there is a sequence of projections $p_{i}:
T_{u_0}P^{i+1} \rightarrow T_{u_0}P^{i}$, with 
 $p_i ( HT_{u_0}P^{i+1}) = HT_{u_0} P^i$ and $p_i ( VT_{u_0}P^{i+1}) = VT_{u_0}
P^i$, where $P^{n+1}=GL(M)$.

    \item[(C2)] The projections
$p_{n} \circ p_{n-1} \circ \ldots p_{i}[\delta  Ad (\xi) X_j (u_0)]$ are in the
image $ Im ( \delta g^i(M))$, for
all vector
fields
$X_j$ , 
$j=0,1,\ldots ,k$ of the sde $(\ref{eq1})$ and for any automorphism $\xi$ of
$P^{i}$.

%Let $P^j$ a $G$-structure associated to tensor fields
%    $k^j$, with $P^1 \subset ... \subset P^{n-1} \subset P^n
%    =GL(M)$, and still,
%    $$ T_{u_0}GL(M)\cong T_{u_0}P^{n-1}\oplus u_0 \qtiu^{n-1}, $$
%    $$ T_{u_0}GL(M)\cong T_{u_0}P^{n-2}\oplus u_0 \qtiu^{n-2}, $$
%    $$\vdots $$
%    $$ T_{u_0}GL(M)\cong T_{u_0}P^{1}\oplus u_0 \qtiu^{1}, $$
%    with $\qtiu^{n-1} \subset \qtiu^{n-2} \subset... \subset \qtiu^1 \subset
%    \gtiu$.
\end{description}

\begin{corollary}[Cascade decomposition] \label{cascata}
Under conditions (C1) and (C2) above, we have the following
decomposition 
of the stochastic flow:
$$
 \varphi_t = \xi_t^1 \circ \xi_t^2 \circ \ldots  \xi_t^n \circ \rho_t,
$$
where for $1\leq i \leq
n$,  $\xi_t^i \in G^k(M)$, the partial compositions  
$(\xi_t^1 \circ ...  \xi_t^i)$ are diffusions in $G^i(M)$ , 
$(\xi_t^{i+1} \circ ...  \xi_t^n \circ \rho_t)(x_0)=x_0$ and the derivative
$(\xi_t^{i+1} \circ ...  \xi_t^n \circ \rho_t)_{*} (u_0)=u_0 q^{(i)}_t$ where
$q^{(i)}_t$ is a process in the
Lie
subgroup $ \langle \exp(\ker p_{n} \circ p_{n-1} \circ \ldots p_{i} )\rangle
\subset \Gl$.

%\forall t \geq 0$. Each component $(\xi_t^1 \circ \xi_t^2 \circ
%... \circ \xi_t^k) \in G^k(M)$, $1\leq k \leq n-1$, $(\xi_t^{k+1}
%\circ ... \circ \xi_t^n \circ \rho_t)(x_0)=x_0$ and $(\xi_t^{k+1}
%\circ ... \circ \xi_t^n \circ \rho_t)_*(u_0)=u_0 q_t$ with $q_t$
%been a process in $\qtiu^k \subset \gtiu$.
\end{corollary}

\proof  For each subbundle $P^i$ of $GL(M)$, consider the projection $p_{n}
\circ p_{n-1} \circ \ldots p_{i}: T_{u_0}GL(M) \rightarrow T_{u_0}P^i $. Condition
(C2) implies hypothesis (H2), hence, by Theorem \ref{geral}
there exists a decomposition $\varphi_t = \tilde{\xi}^{(i)}_t
\rho^{(i)}_t$ such that $\tilde{\xi}_t^{(i)}$ is a diffusion in $G^i(M)$ and 
$\rho^{(i)}_{t*} u_0 = u_0 q^{(i)}_t$. The result follows by taking
$\xi^1_t = \tilde{\xi}^{(1)}_t$ and by induction 
\[
\xi_t^i = (\tilde{\xi}^{(i-1)}_t)^{-1} \circ \tilde{\xi}^{(i)}_t
\] 
for $1<i \leq  n$ and $\rho_t = \rho^{(n)}_t$.

%We decompose
%initially $\varphi_t = \Psi_t \circ \rho_t$, with $\Psi_t$ a diffusion in the
%group of automorphisms of the largest subbundle  $P^n$. The argument now
%follows
%by induction decomposing $\Psi_t$ into a component in $P^{n-1}$, and so on up
%to
%the $P^1$ preserving component.

\eop

%initially the compon The first step is to obtain a decomposition $\phi_t =
%\xi_t
%\circ \rho_t$, with $\xi_t$ been a diffusion in affine
%transformation group and the ``rest" $\rho_t$ had property how
%$\rho_t(x_0)=x_0$ and $\rho_{t*}=Id_{T_{x_0}M}, \ \forall t \geq
%0$. Using the Theorem \ref{geral} to the affine flow $\xi_t$, we
%have that $\xi_t = \xi_t^1 \circ \Delta_t^1$, with $\xi_t^1 \in
%G^1(M)$, $\Delta_t^1 (x_0)=x_0$ and $\Delta_{t*}^1 u_0=u_0q_t$,
%where $q_t \in \qtiu^1$. Furthermore, how
%$\rho_{t*}=Id_{T_{x_0}M}$ we have that $(\Delta_{t}^1 \circ
%\rho_t)_* u_0=u_0q_t$.

%By again Theorem \ref{geral}, $\xi_t= (\xi_t^1 \circ \xi_t^2)
%\circ \Delta_t^2$ with $(\xi_t^1 \circ \xi_t^2) \in G^2(M)$,
%$(\Delta_t^2 \circ \rho_t)x_0=x_0$ e $(\Delta_t^2 \circ \rho_t)_*
%u_0= u_0q_t$, where now $q_t \in \qtiu^2$. Make it successively,
%at $(n-1)$-th step, we will get $\xi_t= (\xi_t^1 \circ \xi_t^2
%\circ...\circ \xi_t^{n-1}) \circ \Delta_t^{n-1}$ with $(\xi_t^1
%\circ \xi_t^2 \circ...\circ \xi_t^{n-1}) \in G^{n-1}(M)$,
%$(\Delta_t^{n-1} \circ \rho_t)x_0=x_0$ and $(\Delta_t^{n-1} \circ
%\rho_t)_* u_0= u_0q_t$, where $q_t \in \qtiu^{n-1}$.

%Finally, more a time by Theorem \ref{geral}, we can obtain $\xi_t=
%(\xi_t^1 \circ \xi_t^2 \circ...\circ \xi_t^{n}) \circ
%\Delta_t^{n}$ with $(\xi_t^1 \circ \xi_t^2 \circ...\circ
%\xi_t^{n}) \in G^{n}(M)$. By uniqueness of decomposition $\phi_t =
%\xi_t \circ \rho_t$, follow that $\Delta_t^n = \rho_t$, and
%therefore $$ \phi_t = \xi_t^1 \circ ... \circ \xi_t^{n-1} \circ
%\xi_t^n \circ \rho_t. $$

A direct  example of cascade decomposition occurs if we put together the
results of Theorems \ref{decomp isom}, \ref{decomp afim} and \ref{pres vol}:
If the geometry of the manifold is simple enough (large groups of
isometries and affine transformations) then the flow decomposes as
\[
 \varphi_t = \xi_t^1 \circ \xi_t^2 \circ \xi_t^3 \circ \rho_t,
\]
where $\xi_t^1$ is a diffusion in the group of isometries of $M$, $\xi_t^1
\circ \xi_t^2$ is a diffusion in the group of volume preserving transformations
of $M$, $\xi_t^1
\circ \xi_t^2 \circ \xi_t^3$ is a diffusion in the group of affine
transformations
of $M$ and 
 $\rho_t (x_0)=x_0$ with 
linearization at $x_0$ given by  $\rho_{t*} = Id$.

%Now, we suppose that a $SO(n,\mathbb{R})$-structure $P^1$ is
%counted in a $Sl(n,\mathbb{R})$-structure $P^2$, i.e.,
%$$ P^1 \subset P^2 \subset GL(M).$$ This way, the metric
%$g=g_{ij}dx_i \otimes dx_j$ on manifold $M$ is associated to
%volume form $\mu_{g}$ by relation
%$$ \mu_{g}=\sqrt{\det(g_{ij})} \ dx_1 \wedge ... \wedge dx_n.$$

%Let $\qtiu^2$ be the Lie subalgebra of diagonal matrix in the form
%$a \cdot Id$, $a \in \mathbb{R}$ and $\qtiu^1$ the Lie subalgebra
%of upper triangular matrix. Therefore, we have
%$$ \qtiu^2 \subset \qtiu^1 \subset \gtiu,$$
%and by Corollary \ref{cascata}, we can obtain the factorization $$
%\phi_t = \xi_t^1 \circ \xi_t^2 \circ \xi_t^3 \circ \rho_t, $$
%where $\xi_t^1$ is a diffusion in isometries transformation group
%of $M$, $(\xi_t^1 \circ \xi_t^2)$ is volume preserving, and
%$(\xi_t^1 \circ \xi_t^2 \circ \xi_t^3)$ is a affine
%transformation. Moreover, $\xi_t^2(x_0)= \xi_t^3(x_0) =
%\rho_t(x_0)=x_0$, $(\xi_t^2 \circ \xi_t^3 \circ \rho_t)_* u_0= u_0
%q_t$ with $q_t \in \mathrm{exp}(\qtiu^1)$ and $( \xi_t^3 \circ
%\rho_t)_* u_0= u_0 q_t$ where $q_t \in \mathrm{exp}(\qtiu^2)$.

\bigskip

Another natural class of examples of cascade decomposition which satisfies
hypotheses (C1) and (C2) is a stochastic flow $\varphi_t$ in $\R^n$ where the
Lie subalgebras $\ptiu_i$ are elements of the form: $\ptiu_1= \ontiu$ and the
sequence is a flag in the subalgebras of the form
\[
 \ptiu_i = \ontiu \oplus \left( \begin{array}{ccc}
                                 B_{k_1} &  & 0 \\
0 & \ddots &  \\
0& 0 & B_{k_s}
                                \end{array} \right)
\]
where $B_{k_j}$ are $k_j \times k_j$ upper triangular
matrices with $\sum_{j=1}^s k_j =n$. One finds many choices of a
chain of subalgebras $\ptiu_1=\ontiu \subset \ptiu_2 \subset \ldots \ptiu_n=
\gtiu$ with natural projections as stated in hypothesis (C2).

\end{document}